\documentclass[10pt]{article}
\usepackage{amssymb,amsmath,textcomp,graphicx,graphics}
\newtheorem{satz}{Theorem}[section]

\newtheorem{bemerk1}{Remark}[section]

\newcommand{\iR}{\mathbb{R}}

\newcommand{\iC}{\mathbb{C}}

\newcommand{\oH}{\hspace*{0.39em}\raisebox{0.6ex}{\textdegree}\hspace{-0.72em}H}

\begin{document}
\begin{center}
{\bf\Large Global strong solvability of a quasilinear subdiffusion
problem}
\end{center}
\vspace{0.7em}
\begin{center}
Rico Zacher
\end{center}
{\footnotesize \noindent {\bf address:} Martin-Luther University
Halle-Wittenberg, Institute of Mathematics, Theodor-Lieser-Strasse
5, 06120 Halle, Germany, E-mail:
rico.zacher@mathematik.uni-halle.de} \vspace{0.7em}
 \vspace{0.7em}
\begin{abstract}
We prove the global strong solvability of a quasilinear
initial-boundary value problem with fractional time derivative of
order less than one. Such problems arise in mathematical physics in
the context of anomalous diffusion and the modelling of dynamic
processes in materials with memory. The proof relies heavily on a
regularity result about the interior H\"older continuity of weak
solutions to time fractional diffusion equations, which has been
proved recently by the author.  We further establish a basic $L_2$
decay estimate for the special case with vanishing external source
term and homogeneous Dirichlet boundary condition.
\end{abstract}
\vspace{0.7em}
\begin{center}
\end{center}

\section{Introduction and main results}
Let $T>0$, $N\ge 2$, and $\Omega\subset \iR^N$ be a bounded domain.
The main purpose of this paper is to prove the global strong
solvability of the time fractional quasilinear problem
\begin{align}
\partial_t^\alpha(u-u_0)-D_i\big(a_{ij}(u)D_j u\big) & = f,\;t\in
(0,T),\,x\in \Omega,\nonumber\\
u & = g,\;t\in (0,T),\,x\in \Gamma, \label{quasiprob}\\
u|_{t=0} & = u_0,\;x\in \Omega,\nonumber
\end{align}
where we use the sum convention. Here $\Gamma=\partial \Omega$, $Du$
denotes the gradient of $u$ with respect to the spatial variables
and $\partial_t^\alpha$ stands for the Riemann-Liouville fractional
derivation operator with respect to time of order $\alpha\in (0,1)$;
it is defined by
\[
\partial_t^\alpha v(t,x)=\partial_t \int_0^t
g_{1-\alpha}(t-\tau)v(\tau,x)\,d\tau,\quad t>0,\,x\in \Omega,
\]
where $g_\beta$ denotes the Riemann-Liouville kernel
\[
g_\beta(t)=\,\frac{t^{\beta-1}}{\Gamma(\beta)}\,,\quad
t>0,\;\beta>0.
\]
The functions $f=f(t,x)$, $g=g(t,x)$, and $u_0=u_0(x)$ are given
data.

During the last decade there has been an increasing interest in time
fractional diffusion equations like (\ref{quasiprob}) and special
cases of it. An important application is the modelling of anomalous
diffusion, see e.g.\ the surveys \cite{Metz}, \cite{Metz2}. In this
context, equations of the type (\ref{quasiprob}) are termed {\em
subdiffusion equations} as the time order $\alpha$ is less than one.
While in normal diffusion (described by the heat equation or more
general parabolic equations), the mean squared displacement of a
diffusive particle behaves like const$\cdot t$ for $t\to\infty$, in
the time fractional case this quantity grows as const$\cdot
t^\alpha$, for which there is evidence in a diverse number of
systems, see \cite{Metz} and the references therein. In the case of
equation (\ref{quasiprob}) the diffusion coefficients are allowed to
depend on the unknown $u$.

Another context where equations of the type (\ref{quasiprob}) and
variants of them appear is the modelling of dynamic processes in
materials with {\em memory}. An example is given by the theory of
heat conduction with memory, see \cite{JanI} and the references
therein. Another application is the following special case of a
model for the diffusion of fluids in porous media with memory, which
has been introduced in \cite{CapuFlow}:
\begin{align}
\partial_t^\alpha(q-q_0)-\mbox{div}\,\big(\kappa(q)D q\big) & = f,\;t\in
(0,T),\,x\in \Omega,\nonumber\\
q & = 0,\;t\in (0,T),\,x\in \Gamma, \label{porousmedia}\\
q|_{t=0} & = q_0,\;x\in \Omega.\nonumber
\end{align}
Here $\alpha\in (0,1)$, $q=q(t,x)$ denotes the pressure of the
fluid, $\kappa=\kappa(q)$ stands for the permeability of the porous
medium, and $f$ is related to external sources in the equation of
balance of mass. Model (\ref{porousmedia}) is obtained by combining
the latter equation with a modified version of Darcy's law for the
mass flux $J$ which reads
\[
J=-\partial_t^{1-\alpha}\big(\kappa(q)Dq\big),
\]
and by assuming that the (average) mass of the fluid is proportional
to the pressure. We refer to \cite{JakuDiss}, where a more general
model is discussed.

We next describe our main result concerning (\ref{quasiprob}).
Letting $p>N+\frac{2}{\alpha}$ we will assume that
\begin{itemize}
\item [{\bf (Q1)}] $f\in L_p([0,T];L_p(\Omega))$, $g\in
B_{pp}^{\alpha(1-\frac{1}{2p})}([0,T];L_p(\Gamma))\cap
L_p([0,T];B_{pp}^{2-\frac{1}{p}}(\Gamma))$, $u_0\in
B_{pp}^{2-\frac{2}{p\alpha}}(\Omega)$, and $u_0=g|_{t=0}$ on
$\Gamma$;
\item [{\bf (Q2)}] $A=(a_{ij})_{i,j=1,\ldots,N}\in C^1(\iR;\mbox{Sym}\{N\})$, and
there exists $\nu>0$ such that $a_{ij}(y)\xi_i\xi_j\ge \nu|\xi|^2$
for all $y\in \iR$ and $\xi\in \iR^N$.
\end{itemize}
Here Sym$\{N\}$ denotes the space of $N$-dimensional real symmetric
matrices. For $s>0$ and $1<p<\infty$ the symbols $H^s_p$ and
$B_{pp}^s$ refer to Bessel potential (Sobolev spaces for integer
$s$) and Sobolev-Slobodeckij spaces, respectively.

The main result concerning the problem (\ref{quasiprob}) reads as
follows.
\begin{satz} \label{quasilinear}
Let $\Omega\subset \iR^N$ ($N\ge 2$) be a bounded domain with
$C^2$-smooth boundary $\Gamma$. Let $\alpha\in (0,1)$, $T>0$ be an
arbitrary number, $p>N+\frac{2}{\alpha}$, and suppose that the
assumptions (Q1) and (Q2) are satisfied. Then the problem
(\ref{quasiprob}) possesses a unique strong solution $u$ in the
class
\[
u\in H^\alpha_p([0,T];L_p(\Omega))\cap L_p([0,T];H^2_p(\Omega)).
\]
\end{satz}
Note that short-time existence of strong or classical solutions to
problems like (\ref{quasiprob}) can be established by means of
maximal regularity for the linearized problem and the contraction
mapping principle. This has been known before, see e.g.\ \cite{CLS},
\cite{ZQ}, and for the semilinear case \cite{CP2}. For results on
maximal regularity for fractional evolution equations we further
refer to \cite{CleGriLon}, \cite{Kochabs}, \cite{JanI}, and
\cite{ZEQ}. The novelty here is that $T>0$ can be given arbitrarily
large without assuming any smallness condition on the data. In some
papers (global) generalized solutions are constructed for
quasilinear subdiffusion problems, see \cite{Grip1} and \cite{Jaku}.
These results are based on the theory of accretive operators.

The crucial step in our
proof of the global existence result is an estimate of the H\"older norm of $u$ on parabolic subdomains $[0,\delta]\times \bar{\Omega}$
which is uniform with respect to $\delta\in (0,T]$. In a very recent work (\cite{ZHoelder}, see also
\cite{ZHabil}), the author was able to prove {\em interior} H\"older regularity for weak solutions of time
fractional diffusion equations of the form
\begin{equation} \label{IntroMProb}
\partial_t^\alpha(u-u_0)-D_i\big(a_{ij}(t,x)D_j u\big)=f,\;\;t\in (0,T),\,x\in \Omega,
\end{equation}
with $\alpha\in (0,1)$ and merely bounded and measurable coefficients, see Theorem \ref{interiorHoelder} below.
Using this result, we derive conditions which are
sufficient for H\"older continuity {\em up to the parabolic boundary}. Here we do not aim at high generality but we
are content with finding some simple
conditions which are also necessary when studying the quasilinear problem (\ref{quasiprob}) in the setting of
maximal $L_p$-regularity. Because
of the latter it is natural and also not so difficult to use the method of maximal $L_p$-regularity to achieve the goal.

In this paper we also prove a basic decay estimate for the $L_2(\Omega)$-norm of the solution $u$ to (\ref{quasiprob}) in
the special case when $f=g=0$. It is shown that for the
global strong solution of (\ref{quasiprob}) we have in this case
\begin{equation} \label{decest}
|u(t,\cdot)|_{L_2(\Omega)}^2\le
\,\frac{c|u_0|_{L_2(\Omega)}^2}{1+\mu t^\alpha},\quad t\ge 0,
\end{equation}
where $c=c(\alpha)$ and $\mu=\mu(\nu,N,\Omega)$ are positive
constants. A polynomial decay estimate has already been known in the
linear case for problems of the type
\begin{align}
\partial_t^\alpha(u-u_0)-D_i\big(a_{ij}(x)D_j u\big) & = 0,\;t\in
(0,T),\,x\in \Omega,\nonumber\\
u & = 0,\;t\in (0,T),\,x\in \Gamma, \label{decayprob}\\
u|_{t=0} & = u_0,\;x\in \Omega,\nonumber
\end{align}
see \cite[Cor.\ 4.1]{Naka}; for the special case
$a_{ij}=\delta_{ij}$ we also refer to \cite[p.\ 11--13]{Meer}. By
means of an eigenfunction expansion for $u$, there it is shown that
in general $|u(t,\cdot)|_{L_2(\Omega)}$ decays as $1/t^\alpha$ as
$t\to \infty$, which is optimal w.r.t.\ the exponent and stronger
than (\ref{decest}), where we have this behaviour for the square of
$|u(t,\cdot)|_{L_2(\Omega)}$. In particular we do not have
exponential decay as in the case $\alpha=1$.

The smaller exponent ($\alpha/2$ instead of $\alpha$) in the
quasilinear case is due to the different method. Note that the
aforementioned method is no longer applicable in the quasilinear
case. Our proof of (\ref{decest}) is based on energy estimates; it
makes use of the fundamental convexity inequality (\ref{ident2}) for
the Riemann-Liouville fractional derivative, see Theorem
\ref{fundament} below. Our method applies to a more general class of
quasilinear problems, including e.g.\ the time fractional
$p$-Laplace equation, and it extends to the (more natural) weak
setting. This will be elaborated in a forth-coming paper.

The paper is organized as follows. In Section 2 we collect some
known results on a priori estimates for linear time fractional
diffusion equations and recall the fundamental convexity identity
(\ref{ident2}) for the fractional derivative. Section 3 is devoted
to H\"older regularity up to $t=0$ for weak solutions of
(\ref{IntroMProb}). Regularity up to the full parabolic boundary is
then established in Section 4. Using these estimates, the global
existence result, Theorem \ref{quasilinear}, is proved in Section 5.
Finally, we derive the described decay estimate in Section 6.
\section{Preliminaries}
We first fix some notation. For $T>0$ and a bounded domain
$\Omega\subset \iR^N$ with boundary $\Gamma$ we put
$\Omega_T=(0,T)\times \Omega$ and $\Gamma_T=(0,T)\times \Gamma$. The
Lebesgue measure in $\iR^N$ will be denoted by $\lambda_N$.

The boundary $\Gamma$ is said to satisfy the property of {\em
positive geometric density}, if there exist $\beta\in (0,1)$ and
$\rho_0>0$ such that for any $x_0\in \Gamma$, any ball $B(x_0,\rho)$
with $\rho\le \rho_0$ we have that $\lambda_N(\Omega\cap
B(x_0,\rho))\le \beta\lambda_N(B(x_0,\rho))$, cf.\ e.g.\
\cite[Section I.1]{DB}.

By $y_+:=\max\{y,0\}$ we denote the positive part of $y\in \iR$.

In the following we collect some known results from the linear
theory which are basic to the investigation of (\ref{quasiprob}).

Let $T>0$ and $\Omega$ be a bounded domain in $\iR^N$ with $N\ge 2$.
We consider the linear time fractional diffusion equation
\begin{equation} \label{MProb}
\partial_t^\alpha(u-u_0)-D_i\big(a_{ij}(t,x)D_j u\big)=f,\;\;t\in (0,T),\,x\in \Omega.
\end{equation}
We assume that
\begin{itemize}
\item [{\bf (H1)}] $A\in L_\infty(\Omega_T;\iR^{N\times
N})$, and
\[
\sum_{i,j=1}^N|a_{ij}(t,x)|^2\le \Lambda^2,\quad \mbox{for
a.a.}\;(t,x)\in \Omega_T.
\]
\item [{\bf (H2)}] There exists $\nu>0$ such that
\[
a_{ij}(t,x)\xi_i\xi_j\ge \nu|\xi|^2,\quad\mbox{for a.a.}\;
(t,x)\in\Omega_T,\; \mbox{and all}\;\xi\in \iR^N.
\]
\item [{\bf (H3)}] $u_0\in L_\infty(\Omega)$; $f\in
L_r([0,T];L_q(\Omega))$, where $r,q\ge 1$ fulfill
\end{itemize}
\[
\frac{1}{\alpha r}\,+\,\frac{N}{2q}\,=1-\kappa,
\]
and
\[
r\in \Big[\,\frac{1}{\alpha(1-\kappa)}\,,\infty\Big],\; q\in
\Big[\,\frac{N}{2(1-\kappa)}\,,\infty\Big],\;\kappa\in (0,1).
\]
Following \cite{ZHoelder} and \cite{ZHabil} we say that a function
$u$ is a {\em weak solution} of (\ref{MProb}) in $\Omega_T$, if $u$
belongs to the space
\begin{align*}
{\cal S}_\alpha:=\{&\,v\in L_{2/(1-\alpha),w}([0,T];L_2(\Omega))\cap
L_2([0,T];H^1_2(\Omega))\;
\mbox{such that}\;\\
&\;\;g_{1-\alpha}\ast v\in C([0,T];L_2(\Omega)),
\;\mbox{and}\;(g_{1-\alpha}\ast v)|_{t=0}=0\},
\end{align*}
and for any test function
\[
\eta\in \oH^{1,1}_2(\Omega_T):=H^1_2([0,T];L_2(\Omega))\cap
L_2([0,T];\oH^1_2(\Omega))
 \quad\;\;
\Big(\oH^1_2(\Omega):=\overline{C_0^\infty(\Omega)}\,{}^{H^1_2(\Omega)}\Big)
\]
with $\eta|_{t=T}=0$ there holds
\begin{equation} \label{BWF}
\int_{0}^{T} \int_\Omega \Big(-\eta_t \big[g_{1-\alpha}\ast
(u-u_0)\big]+ a_{ij}D_j u D_i \eta\Big)\,dxdt = \int_{0}^{T}
\int_\Omega f\eta \,dxdt.
\end{equation}
Here $L_{p,\,w}$ stands for the weak $L_p$ space and $f_1\ast f_2$
means the convolution on the positive halfline with respect to time,
that is $(f_1\ast f_2)(t)=\int_0^t f_1(t-\tau)f_2(\tau)\,d\tau$,
$t\ge 0$.

Existence of weak solutions of (\ref{MProb}) in the class ${\cal
S}_\alpha$ has been shown in \cite{ZWH}. For example, assuming (H1),
(H2), $u_0\in L_2(\Omega)$, and $f\in L_2(\Omega_T)$, the
corresponding Dirichlet problem has a unique solution in ${\cal
S}_\alpha$. Global boundedness of weak solutions has been obtained
in \cite{Za}, the result can be stated as follows, cf.
\cite[Corollary 3.1]{Za}.
\begin{satz} \label{globalbddness}
Let $\alpha\in (0,1)$, $T>0$ and $\Omega$ be a bounded domain in
$\iR^N$ ($N\ge 2$) with $\Gamma$ satisfying the property of positive
density. Let further the assumptions (H1)-(H3) be satisfied. Suppose
that $u\in {\cal S}_\alpha$ is a weak solution of (\ref{MProb}) such
that $|u|\le K$ a.e.\ on $\Gamma_T$ (in the sense that $(u-K)_+,
(-u-K)_+\in L_2([0,T];\oH^1_2(\Omega))$) for some
$K\ge|u_0|_{L_\infty(\Omega)}$. Then $u$ is essentially bounded in
$\Omega_T$ and
\[
|u|_{L_\infty(\Omega_T)}\le C(1+K),
\]
where the constant $C=C(\alpha,r,q,T,N,\nu,\Omega,|f|_{L_r(L_q)})$.
\end{satz}
An {\em interior} H\"older estimate for bounded weak solutions of
(\ref{MProb}) has been proved recently in \cite{ZHoelder}, see also
\cite{ZHabil}. For $\beta_1,\beta_2\in (0,1)$ and $Q\subset
\Omega_T$ we put
\[
[u]_{C^{\beta_1,\beta_2}(Q)}:=\sup_{(t,x),(s,y)\in
{Q},\,(t,x)\neq(s,y)}\Big\{
\frac{|u(t,x)-u(s,y)|}{|t-s|^{\beta_1}+|x-y|^{\beta_2}} \Big\}.
\]
Then the interior regularity result reads as follows, cf.
\cite[Theorem 1.1]{ZHoelder}.
\begin{satz} \label{interiorHoelder}
Let $\alpha\in (0,1)$, $T>0$ and $\Omega$ be a bounded domain in
$\iR^N$ ($N\ge 2$). Let the assumptions (H1)-(H3) be satisfied and
suppose that $u\in {\cal S}_\alpha$ is a bounded weak solution of
(\ref{MProb}) in $\Omega_T$. Then there holds for any
$Q\subset\Omega_T$ separated from $\Gamma_T$ by a positive distance
$d$,
\[
[u]_{C^{\frac{\alpha\epsilon}{2},\epsilon}(\bar{Q})}\le
C\Big(|u|_{L_\infty(\Omega_T)}+|u_0|_{L_\infty(\Omega)}+|f|_{L_r([0,T];L_q(\Omega))}\Big)
\]
with positive constants
$\epsilon=\epsilon(\Lambda,\nu,\alpha,r,q,N,\mbox{diam}\,\Omega,\inf_{(\tau,z)\in
Q}\tau)$ and $C=C(\Lambda,\nu,\alpha,r,q,N$,
$\mbox{diam}\,\Omega,\lambda_{N+1}(Q),d)$.
\end{satz}
\begin{bemerk1} \label{fsumme}
{\em The statement of Theorem \ref{interiorHoelder} can be extended
to the case where the right-hand side of equation (\ref{MProb}) has
the form
\[
\sum_{k=1}^{k_f}f_k-\sum_{k=1}^{k_g}D_i g_k^i,
\]
with $f_k\in L_{r_k}([0,T];L_{q_k}(\Omega))$, $k=1,\ldots,k_f$,
$\sum_{i=1}^N (g_k^i)^2\in L_{r^{(k)}}([0,T];L_{q^{(k)}}(\Omega))$,
$k=1,\ldots,k_g$, and all pairs of exponents $(r_k,q_k)$ and
$(r^{(k)},q^{(k)})$, respectively, are subject to the condition in
(H3). This follows from \cite[Remark 6.1]{ZHoelder}, see also
\cite{ZHabil}.}
\end{bemerk1}
Next we are concerned with maximal $L_p$-regularity for the
corresponding problem in non-divergence form,
\begin{align}
\partial_t^\alpha\big(u-u_0\big)-a_{ij}(t,x)D_i D_j u& =f,\;\;t\in (0,T),\,x\in
\Omega,\nonumber\\
u &=g,\;\;t\in (0,T),\,x\in \Gamma,\label{linearMR}\\
u|_{t=0} &=u_0,\;\;x\in \Omega,\nonumber
\end{align}
where we again use the sum convention. The following result is a
special case of \cite[Theorem 3.4]{ZQ} on linear initial-boundary
value problems in the context of parabolic Volterra equations.
\begin{satz} \label{linearmaxreg}
Let $T>0$ and $\Omega\subset \iR^N$ be a bounded domain with
$C^2$-boundary $\Gamma$, and $N\ge 2$. Let $\alpha\in (0,1)$ and
$p>\frac{1}{\alpha}+\frac{N}{2}$. Suppose that
$A=(a_{ij})_{i,j=1,\ldots,N}\in C([0,T]\times
\overline{\Omega};\mbox{{\em Sym}}\{N\})$, and there exists $\nu>0$
such that $a_{ij}(t,x)\xi_i\xi_j\ge \nu|\xi|^2$ for all $(t,x)\in
[0,T]\times \overline{\Omega}$ and $\xi\in \iR^N$. Then the problem
(\ref{linearMR}) has a unique solution $u$ in the class
\[
Z:=H^\alpha_p([0,T];L_p(\Omega))\cap
L_p([0,T];H^2_p(\Omega))\hookrightarrow C([0,T]\times
\overline{\Omega})
\]
if and only if the following conditions are satisfied.
\begin{itemize}
\item[(i)] $f\in L_p([0,T];L_p(\Omega))$,
$g\in Y_D:=B_{pp}^{\alpha(1-\frac{1}{2p})}([0,T];L_p(\Gamma))\cap
L_p([0,T];B_{pp}^{2-\frac{1}{p}}(\Gamma))$, and $u_0\in Y_\gamma:=
B_{pp}^{2-\frac{2}{p\alpha}}(\Omega)$;
\item[(ii)] $u_0=g|_{t=0}$ on $\Gamma$.
\end{itemize}
In this case one has an estimate of the form
\[
|u|_Z\le C\big(|f|_{L_p(\Omega_T)}+|g|_{Y_D}+|u_0|_{Y_\gamma}\big),
\]
where $C$ only depends on $\alpha,p,N,T,\Omega,A$.
\end{satz}
We conclude this preliminary part with an important convexity
inequality for the Riemann-Liouville fractional derivation operator,
which will be needed in Section \ref{SectionDecay}.
\begin{satz} \label{fundament}
Let $\alpha\in (0,1)$, $T>0$ and ${\cal H}$ be a real Hilbert space
with inner product $(\cdot|\cdot)_{{\cal H}}$. Suppose that $v\in
L_{2}([0,T];{\cal H})$ and that there exists $x\in{\cal H}$ such
that $v-x\in \mbox{}_0H^\alpha_2([0,T];{\cal H}):=\{g_\alpha\ast w:
w\in L_{2}([0,T];{\cal H})\}$. Then
\begin{align}
\Big(v(t),\frac{d}{dt}\,(g_{1-\alpha}\ast v)(t)\Big)_{{\cal H}}
\ge&\; \frac{1}{2}\,\frac{d}{dt}\,(g_{1-\alpha}\ast |v|_{\cal
H}^2)(t)+\frac{1}{2}\,g_{1-\alpha}(t)|v(t)|_{\cal H}^2, \quad
\mbox{a.a.}\;t\in(0,T).\label{ident2}
\end{align}
\end{satz}
{\em Proof.} This follows from Theorem 2.1, Proposition 2.1, and
Example 2.1 in \cite{VZ}. \hfill $\square$
\section{Regularity up to $t=0$}
The objective of this and the following section is to find
conditions on the data which ensure H\"older continuity up to the
parabolic boundary for weak solutions of the linear time fractional
diffusion equation (\ref{MProb}). As already mentioned in the
introduction we do not aim at great generality but at results which
are sufficient for the quasilinear problem to be studied.

We first discuss regularity up to $t=0$.
\begin{satz} \label{Regupt=0}
Let $\alpha\in (0,1)$, $T>0$ and $\Omega$ be a bounded domain in
$\iR^N$ ($N\ge 2$). Let the assumptions (H1)-(H3) be satisfied. Let
further $\Omega'\subset\Omega$ be an arbitrary subdomain and assume
that
\[
u_0|_{\tilde{\Omega}}\in
B_{pp}^{2-\frac{2}{p\alpha}}(\tilde{\Omega})\quad \mbox{with}\quad
p>\frac{1}{\alpha}+\frac{N}{2},
\]
for some $C^2$-smooth domain $\tilde{\Omega}$ such that
$\Omega'\subset\tilde{\Omega}\subset\Omega$ and $\Omega'$ is
separated from $\partial\tilde{\Omega}$ by a positive distance $d$.
Then, for any bounded weak solution $u$ of (\ref{MProb}) in
$\Omega_T$, there holds
\begin{equation} \label{Regt=0est}
[u]_{C^{\frac{\alpha\epsilon}{2},\epsilon}([0,T]\times
\overline{\Omega'}\,)}\le
C\Big(|u|_{L_\infty(\Omega_T)}+|u_0|_{L_\infty(\Omega)}+|u_0|_{B_{pp}^{2-\frac{2}{p\alpha}}(\tilde{\Omega})}+|f|_{L_r([0,T];L_q(\Omega))}\Big)
\end{equation}
with positive constants
$\epsilon=\epsilon(\Lambda,\nu,\alpha,p,r,q,N,\mbox{diam}\,\Omega)$
and $C=C(\Lambda,\nu,\alpha,p,r,q,N,d,\mbox{diam}\,\Omega,T,$
$\lambda_{N}(\Omega'))$.
\end{satz}
{\em Proof.} The basic idea of the proof is to extend $u$ to
$[-1,T]\times \Omega$ such that $u$ is H\"older continuous on
$[-1,0]\times \overline{\Omega'}$ and to apply Theorem
\ref{interiorHoelder}.

To this purpose we first extend $u_0|_{\tilde{\Omega}}\in
B_{pp}^{2-\frac{2}{p\alpha}}(\tilde{\Omega})$ to a function
$\hat{u}_0\in B_{pp}^{2-\frac{2}{p\alpha}}(\iR^N)$. By \cite[Theorem
3.1]{ZQ}, the problem
\begin{align*}
\partial_t^\alpha\big(w-\hat{u}_0\big)-\Delta w & =0,\; \;t\in
(0,1),\,x\in\iR^N,\\
w|_{t=0} & = \hat{u}_0,\;\;x\in \iR^N,
\end{align*}
possesses a unique solution $w$ in the class
\[
Z:=H^\alpha_p([0,1];L_p(\iR^N))\cap L_p([0,1];H^2_p(\iR^N)),
\]
and one has an estimate of the form
\[
|w|_Z\le C_0|\hat{u}_0|_{B_{pp}^{2-\frac{2}{p\alpha}}(\iR^N)}\le
\tilde{C}_0|u_0|_{B_{pp}^{2-\frac{2}{p\alpha}}(\tilde{\Omega})}.
\]
Note that by the mixed derivative theorem (cf.\ \cite{Sob}),
\[
Z\hookrightarrow
H^{\alpha(1-\varsigma)}_p([0,1];H^{2\varsigma}_p(\iR^N)),\quad
\varsigma\in[0,1],
\]
and thus $Z\hookrightarrow BUC^\delta([0,1]\times \iR^N)$ for some
sufficiently small $\delta\in (0,\alpha/2)$. In fact, the assumption
$p>\frac{1}{\alpha}+\frac{N}{2}$ ensures existence of some
$\varsigma\in (0,1)$ with $\alpha(1-\varsigma)-\frac{1}{p}>\delta$
and $2\varsigma-\frac{N}{p}>\delta$.

Multiplying $w$ by a suitable smooth cut-off function $\varphi(t)$
we can construct a function $\hat{w}\in Z$ with $\hat{w}|_{t=0}=
\hat{u}_0$ and $\hat{w}|_{t=1}=0$. We then extend $u$ to
$[-1,T]\times \Omega$ by setting $u(t,x)=\hat{w}(-t,x)$ for $t\in
[-1,0)$ and $x\in\Omega$.

Next, we shift the time by setting $\tau=t+1$. Put
$\hat{u}(\tau,x)=u(\tau-1,x)$, $\tau\in (0,T+1)$, $x\in
\tilde{\Omega}$. Define further
\[
g:=\partial_\tau^\alpha \hat{u}-\Delta \hat{u},\quad
\tau\in(0,1),\,x\in \tilde{\Omega}.
\]
Then $g\in L_p([0,1]\times \tilde{\Omega})$, since
$\hat{u}|_{\tau\in(0,1)}\in
H^\alpha_p([0,1];L_p(\tilde{\Omega}))\cap
L_p([0,1];H^2_p(\tilde{\Omega}))$ and $\hat{u}|_{\tau=0}=0$.
Furthermore we have for any test function $\eta\in
\oH^{1,1}_2([0,T+1]\times \tilde{\Omega})$,
\begin{align} \label{schwach1}
\int_0^1\int_{\tilde{\Omega}}\Big(-\eta_\tau(g_{1-\alpha}\ast
\hat{u})+D_j\hat{u}
D_i\eta\Big)\,dx\,d\tau=\int_0^1\int_{\tilde{\Omega}}g\eta\,dx\,d\tau
-\int_{\tilde{\Omega}}\eta(g_{1-\alpha}\ast
\hat{u})\,dx\Big|_{\tau=1}.
\end{align}

On the other hand, we have for a.a.\ $(\tau,x)\in (1,T+1)\times
\tilde{\Omega}$,
\begin{align*}
(g_{1-\alpha}\ast\hat{u})(\tau,x) &
=(g_{1-\alpha}\ast{u})(\tau-1,x)+
\int_0^1 g_{1-\alpha}(\tau-\sigma)\hat{u}(\sigma,x)\,d\sigma\\
& =
\big(g_{1-\alpha}\ast({u}-u_0)\big)(\tau-1,x)+g_{2-\alpha}(\tau)u_0(x)\\
&\quad\quad +\int_0^1
g_{1-\alpha}(\tau-\sigma)\big(\hat{u}(\sigma,x)-u_0(x)\big)\,d\sigma.
\end{align*}
Set
\[
h(\tau,x)=g_{1-\alpha}(\tau)u_0(x)+\int_0^1
\dot{g}_{1-\alpha}(\tau-\sigma)\big(\hat{u}(\sigma,x)-u_0(x)\big)\,d\sigma
=:h_1(\tau,x)+h_2(\tau,x),
\]
$\hat{a}_{ij}(\tau,x)=a_{ij}(\tau-1,x)$, and
$\hat{f}(\tau,x)=f(\tau-1,x)$ for $(\tau,x)\in (1,T+1)\times
\tilde{\Omega}$. Since $u$ is a weak solution of (\ref{MProb}) in
$\Omega_T$, we thus obtain after a short computation that for any
$\eta\in \oH^{1,1}_2([0,T+1]\times \tilde{\Omega})$ with
$\eta|_{\tau=T+1}=0$
\begin{align}
\int_1^{T+1}\int_{\tilde{\Omega}}\Big(-&\eta_\tau(g_{1-\alpha}\ast
\hat{u})+\hat{a}_{ij}D_j\hat{u}
 D_i\eta\Big)\,dx\,d\tau=\nonumber\\
& \int_1^{T+1}\int_{\tilde{\Omega}}(\hat{f}+h)\eta\,dx\,d\tau
+\int_{\tilde{\Omega}}\eta(g_{1-\alpha}\ast
\hat{u})\,dx\Big|_{\tau=1}. \label{schwach2}
\end{align}

Adding (\ref{schwach1}) and (\ref{schwach2}) shows that $\hat{u}$ is
a weak solution of
\[
\partial_\tau^\alpha\hat{u}-D_i(b_{ij}D_j \hat{u})=
\tilde{f},\;\tau\in (0,T+1),\,x\in \tilde{\Omega},
\]
where
\[
b_{ij}(\tau,x)=\chi_{[0,1]}(\tau)+\chi_{(1,T+1]}(\tau)\hat{a}_{ij}(\tau,x)
\]
and
\[
\tilde{f}(\tau,x)=\chi_{[0,1]}(\tau)g(\tau,x)+\chi_{(1,T+1]}(\tau)
\big(\hat{f}+h\big)(\tau,x).
\]

Evidently, $\chi_{[0,1]}(\tau)g\in L_p([0,T+1]\times
\tilde{\Omega})$ and $\chi_{(1,T+1]}(\tau)\hat{f}\in
L_r([0,T+1];L_q(\tilde{\Omega}))$. Concerning the $h$-term we
clearly have $\chi_{(1,T+1]}(\tau)h_1\in L_\infty([0,T+1]\times
\tilde{\Omega})$. To estimate $\chi_{(1,T+1]}(\tau)h_2$, we employ
the H\"older estimate
\[
|\hat{u}(\sigma,x)-u_0(x)|=|\hat{u}(\sigma,x)-\hat{u}(1,x)|\le
C_1(1-\sigma)^\delta,\quad \sigma\in [0,1],\,x\in \tilde{\Omega},
\]
which results from the embedding $Z\hookrightarrow
BUC^\delta([0,1]\times \iR^N)$ and the construction of $\hat{u}$. It
follows that for $1<\tau=t+1\le 1+T$ and $x\in \tilde{\Omega}$
\begin{align*}
|h_2(\tau,x)| & \le C_1\int_0^1
[-\dot{g}_{1-\alpha}(\tau-\sigma)](1-\sigma)^\delta\,d\sigma\\
& = \,\frac{\alpha C_1}{\Gamma(1-\alpha)}\,\int_0^1
(t+\sigma)^{-1-\alpha}\sigma^\delta\,d\sigma.
\end{align*}
Assuming that $t=\tau-1\in (0,1)$ we then have
\begin{align*}
|h_2(\tau,x)| & \le \,\frac{\alpha
C_1}{\Gamma(1-\alpha)}\,\Big(\int_0^t
(t+\sigma)^{-1-\alpha}\sigma^\delta\,d\sigma+\int_t^1
(t+\sigma)^{-1-\alpha}\sigma^\delta\,d\sigma\Big)\\
 & \le \,\frac{\alpha C_1}{\Gamma(1-\alpha)}\,\Big(\int_0^t
(t+\sigma)^{-1-\alpha}t^\delta\,d\sigma+\int_t^1
\sigma^{-1-\alpha+\delta}\,d\sigma\Big)\\
& \le \,\frac{\alpha C_1}{\Gamma(1-\alpha)}\,t^{-\alpha+\delta}
\Big(\frac{1}{\alpha}+\frac{1}{\alpha-\delta}\Big)\\
& \le 3C_1(\tau-1)^\delta g_{1-\alpha}(\tau-1).
\end{align*}
This shows that $\chi_{(1,T+1]}(\tau)h_2\in
L_{r_0}([0,T+1];L_\infty(\tilde{\Omega}))$ for all $1\le
r_0<\frac{1}{\alpha-\delta}$. In particular we find some
$\hat{r}>\frac{1}{\alpha}$ such that $\chi_{(1,T+1]}(\tau)h_2\in
L_{\hat{r}}([0,T+1];L_\infty(\tilde{\Omega}))$.

All in all we see that $\tilde{f}$ is of the form
$\tilde{f}=\sum_{i=1}^4 \tilde{f}_i$, where $\tilde{f}_i\in
L_{r_i}([0,T+1];L_{q_i}(\tilde{\Omega}))$ with
\[
\frac{1}{\alpha r_i}+\frac{N}{2q_i}<1,\quad i=1,2,3,4.
\]
Hence Theorem \ref{interiorHoelder} and Remark \ref{fsumme} imply
that $\hat{u}$ is H\"older continuous in $[1/2,T+1]\times
\overline{\Omega'}$. This in turn yields H\"older continuity of $u$
in $[0,T]\times \overline{\Omega'}$, and it is not difficult to see
that $u$ is subject to the estimate (\ref{Regt=0est}). \hfill
$\square$
\begin{bemerk1} \label{verall}
{\em It follows from Remark \ref{fsumme} and the proof above, that
Theorem \ref{Regupt=0} can be generalized to the case where the
right-hand side of equation (\ref{MProb}) has the form
\[
\sum_{k=1}^{k_f}f_k-\sum_{k=1}^{k_g}D_i g_k^i,
\]
with $f_k$ and $g_k^i$ as in Remark \ref{fsumme}. }
\end{bemerk1}
\section{Regularity up to the parabolic boundary}
The following result gives conditions on the data which are
sufficient for H\"older continuity on $[0,T]\times
\overline{\Omega}$.
\begin{satz} \label{Regparaboundary}
Let $\alpha\in (0,1)$, $T>0$, $N\ge 2$, and $\Omega\subset \iR^N$ be
a bounded domain with $C^2$-smooth boundary $\Gamma$. Let the
assumptions (H1)-(H3) be satisfied. Suppose further that
\[
u_0\in B_{pp}^{2-\frac{2}{p\alpha}}(\Omega),\quad g\in Y_D:=
B_{pp}^{\alpha(1-\frac{1}{2p})}([0,T];L_p(\Gamma))\cap
L_p([0,T];B_{pp}^{2-\frac{1}{p}}(\Gamma))
\]
with $p>\frac{1}{\alpha}+\frac{N}{2}$, and that the compatibility
condition
\[
u_0=g|_{t=0}\quad \mbox{on}\;\,\,\Gamma
\]
is satisfied. Then for any bounded weak solution $u$ of
(\ref{MProb}) in $\Omega_T$ such that $u=g$ a.e. on $(0,T)\times
\Gamma$, there holds
\begin{equation} \label{Regparabound}
[u]_{C^{\frac{\alpha\epsilon}{2},\epsilon}([0,T]\times
\overline{\Omega}\,)}\le
C\Big(|u|_{L_\infty(\Omega_T)}+|u_0|_{B_{pp}^{2-\frac{2}{p\alpha}}({\Omega})}+|f|_{L_r([0,T];L_q(\Omega))}+|g|_{Y_D}\Big)
\end{equation}
with positive constants
$\epsilon=\epsilon(\Lambda,\nu,\alpha,p,r,q,N,\Omega)$ and
$C=C(\Lambda,\nu,\alpha,p,r,q,N,\Omega,T)$.
\end{satz}
{\em Proof.} By Theorem \ref{linearmaxreg}, the problem
\begin{align*}
\partial_t^\alpha(v-u_0)-\Delta v & =0,\;t\in (0,T),\,x\in \Omega\\
v & =g,\;t\in (0,T),\,x\in \Gamma,\\
v|_{t=0} & = u_0,\;x\in \Omega,
\end{align*}
admits a unique strong solution $v$ in the class
\[
v\in Z:=H^\alpha_p([0,T];L_p(\Omega))\cap L_p([0,T];H^2_p(\Omega))
\]
and
\[
|v|_Z\le
C_0(|u_0|_{B_{pp}^{2-\frac{2}{p\alpha}}({\Omega})}+|g|_{Y_D}),
\]
where $C_0$ only depends on $\alpha,p,N,T,\Omega$. As in the proof
of Theorem \ref{Regupt=0} we see that $v\in C^\delta([0,T]\times
\overline{\Omega})$ for some $\delta>0$. Furthermore, the mixed
derivative theorem implies that
\[
D_i v\in H^{\frac{\alpha}{2}}_p([0,T];L_p(\Omega))\cap
L_p([0,T];H^1_p(\Omega))\hookrightarrow H^{\frac{\alpha
\varsigma}{2}}_p([0,T];H_p^{1-\varsigma}(\Omega))
\]
for all $\varsigma\in [0,1]$. Without restriction of generality we
may assume that
$p\in(\frac{1}{\alpha}+\frac{N}{2},\frac{2}{\alpha}+N)$. With
\[
\tilde{p}:=\frac{\frac{1}{\alpha}+\frac{N}{2}}{\frac{2}{\alpha
p}+\frac{N}{p}-1}\,>\,\frac{1}{\alpha}+\frac{N}{2}
\]
and $\varsigma:=\frac{2}{\alpha p}-\frac{1}{\alpha
\tilde{p}}\in(0,1)$ we then have $H^{\frac{\alpha
\varsigma}{2}}_p([0,T];H_p^{1-\varsigma}(\Omega))\hookrightarrow
L_{2\tilde{p}}(\Omega_T)$, which shows that $|D_i v|^2\in
L_{\tilde{p}}(\Omega_T)$ with
$\frac{1}{\alpha\tilde{p}}+\frac{N}{2\tilde{p}}<1$.

Setting $w=u-v$, $w$ is a bounded weak solution of
\[
\partial_t^\alpha w-D_i(a_{ij}D_j w)=f+D_i(a_{ij}D_j v)-\Delta v,\;t\in (0,T),\,x\in
\Omega,
\]
and $w=0$ a.e. on $(0,T)\times \Gamma$.

Next, let $\Omega_0$ be an arbitrary bounded domain containing
$\Omega$ such that dist$(\Omega,\partial\Omega_0)>0$. We extend
$w,f,a_{ij}$ and $\varphi_i:=D_i v$ to $[0,T]\times \Omega_0$ by
setting $w,f,\varphi_i=0$ and $a_{ij}=\delta_{ij}$ on $[0,T]\times
(\Omega_0\setminus \Omega)$. Then $w$ solves
\[
\partial_t^\alpha w-D_i(a_{ij}D_j w)=f+D_i(a_{ij}\varphi_j-\varphi_i),\;t\in (0,T),\,x\in
\Omega_0,
\]
in the weak sense, and thus Theorem \ref{Regupt=0} and Remark
\ref{verall} imply that $w$ is H\"older continuous on $[0,T]\times
\overline{\Omega}$. Since $u=v+w$, the assertion of Theorem
\ref{Regparaboundary} follows. \hfill $\square$
\section{Proof of the global solvability theorem}
The proof of Theorem \ref{quasilinear} is divided into three parts,
devoted respectively to local well-posedness, existence of a
maximally defined solution, and to a priori estimates which lead to
global existence.

Recall that the data belong to the following regularity classes:
\begin{align*}
f\in X^T:=L_p([0,T];L_p(\Omega)),\;\; u_0\in Y_\gamma:=
B_{pp}^{2-\frac{2}{p\alpha}}(\Omega)\\
g\in Y_D^T:= B_{pp}^{\alpha(1-\frac{1}{2p})}([0,T];L_p(\Gamma))\cap
L_p([0,T];B_{pp}^{2-\frac{1}{p}}(\Gamma)).
\end{align*}
We seek a unique solution $u$ of (\ref{quasiprob}) in the space
\[
Z^T:=H^\alpha_p([0,T];L_p(\Omega))\cap L_p([0,T];H^2_p(\Omega)).
\]

{\bf 1. Local well-posedness.} Short-time existence and uniqueness
in the regularity class $Z^\delta$ can be established by means of
the contraction mapping principle and maximal $L_p$-regularity for
an appropriate linearized problem. We proceed similarly as in
\cite{ZQ}, see also \cite{CleLi} and \cite{JanMR}.

We first define a reference function $w\in Z^T$ as the unique
solution of the linear problem
\begin{align*}
\partial_t^\alpha(w-u_0)-a_{ij}(u_0)D_i D_j w & = f+a_{ij}'(u_0)D_i u_0 \,D_j u_0,\,\;t\in
(0,T),\,x\in \Omega,\\
w & = g,\;t\in (0,T),\,x\in \Gamma, \\
w|_{t=0} & = u_0,\;x\in \Omega,
\end{align*}
see Theorem \ref{linearmaxreg}. Note that the condition
$p>N+\frac{2}{\alpha}$ ensures the embedding
\[
u_0\in Y_\gamma= B_{pp}^{2-\frac{2}{p\alpha}}(\Omega)\hookrightarrow
C^1(\overline{\Omega}),
\]
and thus we also have
\[
Z^T\hookrightarrow C([0,T];Y_\gamma)\hookrightarrow
C([0,T];C^1(\overline{\Omega})).
\]

For $\delta\in (0,T]$ and $\rho>0$ let
\[
\Sigma(\delta,\rho)=\{v\in Z^\delta:
v|_{t=0}=u_0,\,|v-w|_{Z^\delta}\le \rho\},
\]
which is a closed subset of $Z^\delta$. By Theorem
\ref{linearmaxreg}, we may define the mapping
$F:\Sigma(\delta,\rho)\rightarrow Z^\delta$ which assigns to $u\in
\Sigma(\delta,\rho)$ the unique solution $v=F(u)$ of the linear
problem
\begin{align}
\partial_t^\alpha(v-u_0)-a_{ij}(u_0)D_i D_j v & = f+h(u,Du,D^2 u),\,\;t\in
(0,\delta),\,x\in \Omega,\nonumber\\
v & = g,\;t\in (0,\delta),\,x\in \Gamma, \label{fixedpoint}\\
v|_{t=0} & = u_0,\;x\in \Omega,\nonumber
\end{align}
where
\[
h(u,Du,D^2 u)=\big(a_{ij}(u)-a_{ij}(u_0)\big)D_i D_j u+a_{ij}'(u)D_i
u\,D_j u.
\]
Observe that every fixed point $u$ of $F$ is a local solution of
(\ref{quasiprob}) and vice versa, at least for some small time
interval $[0,\delta]$.

Since $Z^\delta\hookrightarrow C([0,\delta];C^1(\overline{\Omega}))$
we may set
\[
\mu_w(\delta):=\max\{|w(t,x)-u_0(x)|+|Dw(t,x)-Du_0(x)|:t\in
[0,\delta],\,x\in \overline{\Omega}\}.
\]
Evidently, $\mu_w(\delta)\rightarrow 0$ as $\delta\rightarrow 0$,
due to $w|_{t=0}=u_0$. Letting $u\in \Sigma(\delta,\rho)$ we then
have for all $t\in [0,\delta]$ and $x\in \overline{\Omega}$
\begin{align}
|u(t,x)-u_0(x)|+|Du & (t,x)-Du_0(x)|\le
|u-w|_{C([0,\delta];C^1(\overline{\Omega}))}+\mu_w(\delta)
\nonumber\\
& \le M_0|u-w|_{Z^\delta}+\mu_w(\delta)\le M_0\rho+\mu_w(\delta),
\label{fpa1}
\end{align}
where the embedding constant $M_0>0$ does not depend on $u$ and
$\delta\in (0,T]$; the latter is true since $u-w$ belongs to the
space ${}_0 Z^\delta:=\{\varphi\in Z^\delta: \varphi|_{t=0}=0\}$.
(\ref{fpa1}) yields for any $u\in \Sigma(\delta,\rho)$ the bound
\begin{equation} \label{fpa2}
|u(t,x)-u_0(x)|+|Du(t,x)-Du_0(x)|\le M_0\rho_0+\mu_w(T),\quad t\in
[0,\delta],\,x\in \overline{\Omega},
\end{equation}
where we assume $\rho\in (0,\rho_0]$.

Let now $u\in \Sigma(\delta,\rho)$ and $v=F(u)$. Then $v-w\in {}_0
Z^\delta$ solves the problem
\begin{align}
\partial_t^\alpha(v-w)-a_{ij}(u_0)D_i D_j (v-w) & = h(u,Du,D^2 u)
-a_{ij}'(u_0)D_i u_0 \,D_j u_0,\,\;t\in
(0,\delta),\,x\in \Omega,\nonumber\\
v-w & = 0,\;t\in (0,\delta),\,x\in \Gamma, \nonumber\\
(v-w)|_{t=0} & = 0,\;x\in \Omega.\nonumber
\end{align}
Consequently, it follows from Theorem \ref{linearmaxreg} that for
some constant $M_1>0$ which is independent of $\delta\in (0,T]$
there holds
\begin{align*}
|v-w|_{Z^\delta} & \le M_1 |h(u,Du,D^2 u)-a_{ij}'(u_0)D_i u_0 \,D_j
u_0|_{X^\delta}\\
& \le M_1|\big(a_{ij}(u)-a_{ij}(u_0)\big)D_i D_j
u|_{X^\delta}+M_1|a_{ij}'(u)D_i u\,D_j u-a_{ij}'(u_0)D_i u_0 \,D_j
u_0|_{X^\delta}.
\end{align*}
Using (\ref{fpa1}) and (\ref{fpa2}) we may estimate the first term
as follows.
\begin{align*}
|\big(a_{ij}(u)-a_{ij}(u_0)\big)D_i D_j u|_{X^\delta} & \le
\big(|A(u)-A(w)|_{(L_\infty)^{N^2}}+ |A(w)-A(u_0)|_{(L_\infty)^{N^2}}  \big)\\
& \quad\quad \times \big(|D^2 u-D^2 w|_{(X^\delta)^{N^2}}+|D^2
w|_{(X^\delta)^{N^2}}\big)\\
& \le M_2\big(\rho+\mu_w(\delta)\big)\big(\rho+|D^2
w|_{(X^\delta)^{N^2}}\big),
\end{align*}
where $M_2>0$ does not depend on $\delta$ and $\rho$. Similarly we
obtain
\begin{align*}
|a_{ij}'(u)D_i u\,D_j u-a_{ij}'(u_0)D_i u_0 \,D_j u_0|_{X^\delta}\le
M_3\big(\rho+\mu_w(\delta)\big)\big(\rho+\delta^{\frac{1}{p}}\big),
\end{align*}
with $M_3>0$ being independent of $\delta$ and $\rho$; here the
factor $\delta^{\frac{1}{p}}$ comes from the estimate
$|z|_{X^\delta}\le
\big(\lambda_N(\Omega)\delta\big)^{1/p}|z|_\infty$. We conclude that
\begin{equation} \label{fpa3}
|v-w|_{Z^\delta}\le M\big((\rho+\mu(\delta)\big)^2,
\end{equation}
where $M$ and $\mu(\delta)$ are constants, which do not depend on
$\rho$, $M$ is independent of $\delta$, and $\mu(\delta)$ is
non-decreasing with $\mu(\delta)\rightarrow 0$ as $\delta\rightarrow
0$.

Next let $u_i\in \Sigma(\delta,\rho)$ and $v_i=F(u_i)$, $i=1,2$.
Then $v_1-v_2\in {}_0 Z^\delta$ solves the problem
\begin{align}
\partial_t^\alpha(v_1-v_2)-a_{ij}(u_0)D_i D_j (v-w) & = h(u_1,Du_1,D^2
u_1)-h(u_2,Du_2,D^2 u_2) ,\,\;t\in
(0,\delta),\,x\in \Omega,\nonumber\\
v_1-v_2 & = 0,\;t\in (0,\delta),\,x\in \Gamma, \nonumber\\
(v_1-v_2)|_{t=0} & = 0,\;x\in \Omega,\nonumber
\end{align}
hence
\[
|v_1-v_2|_{Z^\delta}\le M_1|h(u_1,Du_1,D^2 u_1)-h(u_2,Du_2,D^2
u_2)|_{X^\delta}.
\]
Estimating similarly as above we obtain
\begin{equation} \label{fpa4}
|v_1-v_2|_{Z^\delta}\le
M\big((\rho+\mu(\delta)\big)|u_1-u_2|_{Z^\delta},
\end{equation}
where $M$ and $\mu(\delta)$ are like those in (\ref{fpa3}).

Finally, the estimates (\ref{fpa3}) and (\ref{fpa4}) show that for
sufficiently small $\rho$ and $\delta$ the mapping $F$ is a strict
contraction which leaves the set $\Sigma(\delta,\rho)$ invariant.
Local existence and uniqueness of strong solutions to
(\ref{quasiprob}) follows now by the contraction mapping principle.

{\bf 2. The maximally defined solution.} The local solution $u\in
Z^\delta$ obtained in the first part can be continued to some larger
interval $[0,\delta+\delta_1]\subset [0,T]$. In fact, let
$u_\delta:=u|_{t=\delta}\in Y_\gamma$ and define the set
\[
\Sigma(\delta,\delta_1,\rho):=\{v\in
Z^{\delta+\delta_1}:v|_{[0,\delta]}=u,\,|v-w|_{Z^{\delta+\delta_1}}\le
\rho \},
\]
where the reference function $w\in Z^T$ is now defined as the
solution of the linear problem
\begin{align*}
\partial_t^\alpha(w-u_0)-a_{ij}(u_\delta)D_i D_j w & = f+h_1+\chi_{(\delta,T]}(t)a_{ij}'(u_\delta)D_i u_\delta\,D_j u_\delta,\,\;t\in
(0,T),\,x\in \Omega,\\
w & = g,\;t\in (0,T),\,x\in \Gamma, \\
w|_{t=0} & = u_0,\;x\in \Omega,
\end{align*}
with
\[
h_1=\chi_{[0,\delta]}(t)\Big(\big(a_{ij}(u)-a_{ij}(u_\delta)\big)D_i
D_j u+a_{ij}'(u)D_i u\,D_j u\Big).
\]
Observe that $w|_{[0,\delta]}=u$, by uniqueness. So
$\Sigma(\delta,\delta_1,\rho)$ is not empty and it becomes a
complete metric space when endowed with the metric induced by the
norm of $Z^{\delta+\delta_1}$.

Define next the mapping $F:\Sigma(\delta,\delta_1,\rho)\rightarrow
Z^{\delta+\delta_1}$ which assigns to
$\tilde{u}\in\Sigma(\delta,\delta_1,\rho)$ the solution
$v=F(\tilde{u})$ of the linear problem
\begin{align*}
\partial_t^\alpha(v-u_0)-a_{ij}(u_\delta)D_i D_j v & = f+\tilde{h}(\tilde{u},D\tilde{u},D^2\tilde{u})
,\,\;t\in
(0,\delta+\delta_1),\,x\in \Omega,\\
w & = g,\;t\in (0,\delta+\delta_1),\,x\in \Gamma, \\
w|_{t=0} & = u_0,\;x\in \Omega,
\end{align*}
where
\[
\tilde{h}(\tilde{u},D\tilde{u},D^2\tilde{u})=\big(a_{ij}(\tilde{u})-a_{ij}(u_\delta)\big)D_i
D_j \tilde{u}+a_{ij}'(\tilde{u})D_i \tilde{u}\,D_j \tilde{u}.
\]
Since $\tilde{u}|_{[0,\delta]}=u$ we have also $v|_{[0,\delta]}=u$,
by uniqueness.

Observe that $h_1=\tilde{h}(\tilde{u},D\tilde{u},D^2\tilde{u})$ on
$[0,\delta]$ and thus
\[
|v-w|_{Z^{\delta+\delta_1}}\le
M_1|\tilde{h}(\tilde{u},D\tilde{u},D^2\tilde{u})-a_{ij}'(u_\delta)D_i
u_\delta\,D_j u_\delta|_{L_p([\delta,\delta+\delta_1]\times
\Omega)}.
\]
Further,
\[
|F(\tilde{u}_1)-F(\tilde{u}_2)|_{Z^{\delta+\delta_1}}\le
M_1|\tilde{h}(\tilde{u}_1,D\tilde{u}_1,D^2\tilde{u}_1)-\tilde{h}(\tilde{u}_2,D\tilde{u}_2,D^2\tilde{u}_2)|
_{L_p([\delta,\delta+\delta_1]\times \Omega)},
\]
for $\tilde{u}_1,\tilde{u}_2\in \Sigma(\delta,\delta_1,\rho)$.
Therefore we may estimate analogously to the first step to see that
for sufficiently small $\delta_1$ and $\rho$ we have
$F(\Sigma(\delta,\delta_1,\rho))\subset
\Sigma(\delta,\delta_1,\rho)$ and $F$ is a strict contraction. Hence
the contraction mapping principle yields existence of a unique fixed
point of $F$ in $\Sigma(\delta,\delta_1,\rho)$, which is the unique
solution of (\ref{quasiprob}) on $[0,\delta+\delta_1]$.

Repeating this argument we obtain a maximal interval of existence
$[0,T_{max})$ with $T_{max}\le T$, that is $T_{max}$ is the supremum
of all $\tau\in (0,T)$ such that the problem (\ref{quasiprob}) has a
unique solution $u\in Z^{\tau}$.

{\bf 3. A priori bounds and global well-posedness.} In order to
establish global existence we will show that $|u|_{Z^{\tau}}$ stays
bounded as $\tau\nearrow T_{max}$.

Let $\tau\in (0,T_{max})$ and $u\in Z^\tau$ be the unique solution
of (\ref{quasiprob}). Setting $b_{ij}(t,x)=a_{ij}(u(t,x))$, it is
evident that $u$ is a weak solution of
\[
\partial_t^\alpha(u-u_0)-D_i(b_{ij}D_j u)=f,\;t\in (0,\tau),\,x\in \Omega.
\]
Since $Y_\gamma\hookrightarrow C(\overline{\Omega})$ and
$Y_D^\tau\hookrightarrow C([0,\tau]\times \Gamma)$, Theorem
\ref{globalbddness} implies a uniform sup-bound for $|u|$, namely
\[
|u(t,x)|\le C_1,\quad t\in [0,\tau],\,x\in \overline{\Omega},
\]
where the constant $C_1$ depends only on the data
$|f|_{X^T},|g|_\infty,|u_0|_\infty,\Omega,T,\alpha,N$, and $\nu$,
not on $\tau$. It follows then from Theorem \ref{Regparaboundary}
that for some $\varepsilon>0$ we have
\[
|u|_{C^\varepsilon([0,\tau]\times \overline{\Omega})}\le C_2,
\]
where the number $C_2\ge 1$ depends only on
$|f|_{X^T},|g|_{Y_D^T},|u_0|_{Y_\gamma},\Omega,T,\alpha,N$, and
$\nu$, not on $\tau$. In particular, we obtain a uniform H\"older
estimate for the coefficients $b_{ij}$, $i,j=1,\ldots,N$.

The first equation of (\ref{quasiprob}) can be rewritten as
\[
\partial_t^\alpha(u-u_0)-b_{ij}D_i D_j u=f+a_{ij}'(u)D_i u\,D_j u.
\]
By Theorem \ref{linearmaxreg}, the linear problem
\begin{align*}
\partial_t^\alpha(v-u_0)-b_{ij}D_i D_j v & = f,\,\;t\in
(0,\tau),\,x\in \Omega,\\
v & = g,\;t\in (0,\tau),\,x\in \Gamma, \\
v|_{t=0} & = u_0,\;x\in \Omega,
\end{align*}
has a unique solution $v\in Z^\tau$ and there exists a constant
$M_1>0$ independent of $\tau$ such that
\begin{align}
|u-v|_{Z^\tau}& \le M_1|a_{ij}'(u)D_i u\,D_j u|_{X^\tau}\nonumber\\
& \le M_1 \sum_{i,j=1}^N\max_{|y|\le
C_1}|a'_{ij}(y)|\,\big||Du|^2\big|_{X^\tau}. \label{last1}
\end{align}
The assumption on $p$ implies $p>\frac{N}{2}$ and thus
\[
H^2_p(\Omega)\hookrightarrow H^1_{2p}(\Omega)\hookrightarrow
C^{\varepsilon_0}(\overline{\Omega})
\]
for some $\varepsilon_0\in (0,\varepsilon]$. By the
Gagliardo-Nirenberg inequality, there exists then $\theta\in
(0,\frac{1}{2})$ such that
\[
|Du(t,\cdot)|_{L_{2p}(\Omega;\iR^N)}\le
C|u(t,\cdot)|_{H^2_p(\Omega)}^\theta
|u(t,\cdot)|_{C^\varepsilon(\overline{\Omega})}^{1-\theta}\le
CC_2|u(t,\cdot)|_{H^2_p(\Omega)}^\theta,\quad t\in [0,\tau],
\]
and hence by H\"older's and Young's inequality
\begin{align*}
\big||Du|^2\big|_{X^\tau} & \le
\tilde{C}|Du|^2_{L_{2p}([0,\tau]\times \Omega;\iR^N)}\le
C_3|u|_{L_p([0,\tau];H^2_p(\Omega))}^{2\theta}\tau^{\frac{1-2\theta}{p}}\\
& \le C_4|u|_{Z^\tau}^{2\theta}\le
\varepsilon_1|u|_{Z^\tau}+C_5(\varepsilon_1,\theta,C_4),
\end{align*}
for all $\varepsilon_1>0$. This together with (\ref{last1}) yields a
bound for $|u-v|_{Z^\tau}$ which is uniform w.r.t. $\tau$. Since
$|v|_{Z^\tau}$ stays bounded as $\tau\nearrow T_{max}$, it follows
that $|u|_{Z^\tau}$ enjoys the same property. Hence we have global
existence. \hfill $\square$
\section{Decay estimate} \label{SectionDecay}
In this section we prove an $L_2$ decay estimate for solutions of
(\ref{quasiprob}) with $f=0$ and $g=0$, that is we consider
\begin{align}
\partial_t^\alpha(u-u_0)-D_i\big(a_{ij}(u)D_j u\big) & = 0,\;t>0,\,x\in \Omega,\nonumber\\
u & = 0,\;t>0,\,x\in \Gamma, \label{dquasiprob}\\
u|_{t=0} & = u_0,\;x\in \Omega.\nonumber
\end{align}
\begin{satz} Let $\Omega\subset \iR^N$ ($N\ge 2$) be a bounded domain with
$C^2$-smooth boundary $\Gamma$. Let $\alpha\in (0,1)$,
$p>N+\frac{2}{\alpha}$, and $u_0\in
B_{pp}^{2-\frac{2}{p\alpha}}(\Omega)$ such that $u_0=0$ on $\Gamma$.
Assume that condition (Q2) is satisfied. Then for the global strong
solution $u$ of (\ref{dquasiprob}) the function $\{t\mapsto
|u(t,\cdot)|_{L_2(\Omega)}^2\}$ is continuous on $[0,\infty)$ and we
have
\[
|u(t,\cdot)|_{L_2(\Omega)}^2\le
\,\frac{c|u_0|_{L_2(\Omega)}^2}{1+\mu t^\alpha},\quad t\ge 0,
\]
with positive constants $c=c(\alpha)$ and $\mu=\mu(\nu,N,\Omega)$.
\end{satz}
{\em Proof.} Let $T>0$ be arbitrary and $u\in Z^T$ be the solution
of (\ref{dquasiprob}) on $[0,T]$. We multiply the first equation in
(\ref{dquasiprob}) by $u$, integrate over $\Omega$, and integrate by
parts. Using the Dirichlet boundary condition this yields
\[
I_1(t):=\int_\Omega \Big(u\partial_t^\alpha u+a_{ij}(u)D_j u D_i
u\Big)\,dx=\int_\Omega uu_0 g_{1-\alpha}(t)\,dx=:I_2(t),\quad
\mbox{a.a}\;t\in (0,T).
\]
Thanks to (Q2) and Theorem \ref{fundament} with
$\mathcal{H}=L_2(\Omega)$ we have for a.a.\ $t\in (0,T)$
\[
I_1(t)\ge \,\frac{1}{2}\,\partial_t^\alpha
|u(t,\cdot)|_{L_2(\Omega)}^2+\,\frac{1}{2}\,g_{1-\alpha}(t)|u(t,\cdot)|_{L_2(\Omega)}^2+\nu|Du(t,\cdot)|_{L_2(\Omega;\iR^N)}^2.
\]
On the other hand, Young's inequality implies that
\[
I_2(t)\le
\,\frac{1}{2}\,g_{1-\alpha}(t)|u(t,\cdot)|_{L_2(\Omega)}^2+\,\frac{1}{2}\,g_{1-\alpha}(t)|u_0|_{L_2(\Omega)}^2,\quad
\mbox{a.a.}\;t\in (0,T).
\]
Combining these estimates gives
\[
\partial_t^\alpha
|u(t,\cdot)|_{L_2(\Omega)}^2+2\nu|Du(t,\cdot)|_{L_2(\Omega;\iR^N)}^2\le
g_{1-\alpha}(t)|u_0|_{L_2(\Omega)}^2,\quad \mbox{a.a.}\;t\in (0,T),
\]
which in turn, by Poincar\'e's inequality, implies
\begin{equation} \label{diffinequ}
\partial_t^\alpha
|u(t,\cdot)|_{L_2(\Omega)}^2+\mu|u(t,\cdot)|_{L_2(\Omega)}^2\le
g_{1-\alpha}(t)|u_0|_{L_2(\Omega)}^2,\quad \mbox{a.a.}\;t\in (0,T),
\end{equation}
where $\mu=\mu(\nu,N,\Omega)$ is a positive constant. Setting
\[
W(t)=|u(t,\cdot)|_{L_2(\Omega)}^2\quad \mbox{and}\;\;
W_0=W(0)=|u_0|_{L_2(\Omega)}^2
\]
the fractional differential inequality (\ref{diffinequ}) is
equivalent to
\begin{equation} \label{Wineq}
\partial_t^\alpha\big(W-W_0)+\mu W\le 0\quad \mbox{a.e. on} \;(0,T).
\end{equation}

Next, let $V$ denote the solution of the corresponding fractional
differential equation, that is
\begin{equation} \label{Veq}
\partial_t^\alpha\big(V-V_0)+\mu V= 0\quad \mbox{a.e. on}
\;(0,T),\quad V(0)=V_0=W_0.
\end{equation}
By the comparison principle for linear fractional differential
equations (cf. \cite{GLS}), we then have
\[
W(t)\le V(t),\quad t\in [0,T].
\]
The solution of (\ref{Veq}) is given by
\[
V(t)=V_0 E_\alpha(-\mu t^\alpha),\quad t\in [0,T],
\]
where $E_\alpha$ denotes the Mittag-Leffler function defined by
\[
E_\alpha(z)=\sum_{k=0}^\infty \,\frac{z^k}{\Gamma(\alpha k+1)},\quad
z\in \iC,
\]
see \cite[Section 4.1]{KST}. Note that $E_1(z)=e^z$. It is known
that for $\alpha\in (0,1)$ $E_\alpha$ is a completely monotonic
function in $(-\infty,0]$ (see e.g.\ \cite{Poll}) and that there
exists a constant $c>0$ such that
\[
E_\alpha(-x)\le \,\frac{c}{1+x},\quad x\ge 0,
\]
see \cite[Formula (13)]{Kraeh}. It follows that
\[
|u(t,\cdot)|_{L_2(\Omega)}^2\le |u_0|_{L_2(\Omega)}^2E_\alpha(-\mu
t^\alpha)\le \,\frac{c|u_0|_{L_2(\Omega)}^2}{1+\mu t^\alpha},\quad
t\in [0,T].
\]
Since $T>0$ was arbitrary, the theorem is proved. \hfill $\square$



\begin{thebibliography}{99}
{\footnotesize
\bibitem{CapuFlow}
Caputo, M.: Diffusion of fluids in porous media with memory.
Geothermics {\bf 28} (1999), 113-130.
\bibitem{CleLi} Cl\'{e}ment, Ph.; Li, S.: Abstract parabolic
quasilinear evolution equations and applications to a groundwater
problem. Adv. Math. Sci. Appl. {\bf 3} (1994), 17--32.
\bibitem{CleGriLon}
Cl\'{e}ment, Ph.; Gripenberg, G.; Londen, S.-O.: Regularity
properties of solutions of fractional evolution equations. Evolution
equations and their applications in physical and life sciences (Bad
Herrenalb, 1998),  235--246, Lecture Notes in Pure and Appl. Math.,
{\bf 215}, Dekker, New York, 2001.
\bibitem{CLS} Cl\'{e}ment, Ph.; Londen, S.-O.; Simonett, G.:
Quasilinear evolutionary equations and continuous interpolation
spaces. J. Differ. Eq. {\bf 196} (2004), 418--447.
\bibitem{CP2} Cl\'{e}ment, Ph.; Pr\"uss, J.: Global existence for a
semilinear parabolic Volterra equation. Math. Z. {\bf 209} (1992),
17--26.
\bibitem{DB} DiBenedetto, E.: {\em Degenerate parabolic
equations}. Springer, New York, 1993.
\bibitem{Grip1} Gripenberg, G.: Volterra integro-differential
equations with accretive nonlinearity. J. Differ. Eq. {\bf 60}
(1985), 57--79.
\bibitem{GLS} Gripenberg, G.; Londen, S.-O.; Staffans, O.:
{\em Volterra integral and functional equations.} Encyclopedia of
Mathematics and its Applications, {\bf 34}. Cambridge University
Press, Cambridge, 1990.
\bibitem{JakuDiss}
Jakubowski, V.\ G.: {\em Nonlinear elliptic-parabolic
integro-differential equations with $L_1$-data: existence,
uniqueness, asymptotics}. Dissertation, University of Essen, 2001.
\bibitem{Jaku}
Jakubowski, V.\ G.; Wittbold, P.: On a nonlinear elliptic-parabolic
integro-differential equation with $L^1$-data. J. Differential
Equations {\bf 197} (2004), 427--445.
\bibitem{KST} Kilbas, A.\ A.; Srivastava, H.\ M.; Trujillo, J.\ J.:
{\em Theory and applications of fractional differential equations}.
Elsevier, 2006.
\bibitem{Kochabs}
Kochubei, A.\ N.: The Cauchy problem for evolution equations of
fractional order. Differential Equations {\bf 25} (1989), 967--974.
\bibitem{Kraeh} Kr\"ageloh, A.\ M.: Two families of functions related to the
fractional powers of generators of strongly continuous contraction
semigroups. J. Math. Anal. Appl. {\bf 283} (2003), 459-–467.
\bibitem{Meer} Meerschaert, M.\ M., Nane, E., Vellaisamy, P.:
Fractional Cauchy problems on bounded domains. The Annals of
Probability {\bf 37} (2009), 979--1007.
\bibitem{Metz} Metzler, R.; Klafter, J.: The random walk's guide to
anomalous diffusion: a fractional dynamics approach. Phys. Rep. {\bf
339} (2000), 1--77.
\bibitem{Metz2} Metzler, R.; Klafter, J.:
The restaurant at the end of the random walk: recent developments in
the description of anomalous transport by fractional dynamics. J.\
Phys.\ A: Math.\ Gen. {\bf 37} (2004), R161--R208.
\bibitem{Naka} Nakagawa, J., Sakamoto, K., Yamamoto, M.: Overview to
mathematical analysis for fractional diffusion equations -- new
mathematical aspects motivated by industrial collaboration. J.\
Math-for-Industry {\bf 2} (2010A-10), 99--108.
\bibitem{Poll} Pollard, H.: The completely monotonic character of the
Mittag-Leffler function $E_a(-x)$. Bull. Amer. Math. Soc. {\bf 54}
(1948). 1115–-1116.
\bibitem{JanI} Pr\"uss, J.: {\em Evolutionary Integral Equations and
Applications}. Monographs in Mathematics {\bf 87}, Birkh\"auser,
Basel, 1993.
\bibitem{JanMR} Pr\"uss, J.: Maximal regularity for evolution equations in
$L_p$-spaces. Conf. Semin. Mat. Univ. Bari  {\bf 285} (2002), 1--39.
\bibitem{Sob} Sobolevskii, P.E.: Coerciveness inequalities for
abstract parabolic equations. Soviet Math. (Doklady) {\bf 5} (1964),
894--897.
\bibitem{VZ} Vergara, V.; Zacher, R.: Lyapunov functions and
convergence to steady state for differential equations of fractional
order. Math. Z. {\bf 259} (2008), 287--309.
\bibitem{ZHoelder} Zacher, R.: A De Giorgi-Nash type theorem for time fractional diffusion
equations. Submitted to Math.\ Ann.
\bibitem{Za} Zacher, R.: Boundedness of weak solutions to evolutionary partial
integro-differential equations with discontinuous coefficients. J.
Math. Anal. Appl. {\bf 348} (2008), 137--149.
\bibitem{ZHabil}
Zacher, R.: {\em De Giorgi-Nash-Moser estimates for evolutionary
partial integro-differential equations}. Habilitation Thesis. MLU
Halle, 2010. Online available at
http://digital.bibliothek.uni-halle.de/hs/content/titleinfo/825947.
\bibitem{ZEQ} Zacher, R.: Maximal regularity of type $L_p$ for
abstract parabolic Volterra equations. J. Evol. Equ. {\bf 5} (2005),
79--103.
\bibitem{ZQ}
Zacher, R.: Quasilinear parabolic integro-differential equations
with nonlinear boundary conditions.  Differential Integral Equations
{\bf 19} (2006), 1129--1156.
\bibitem{ZWH} Zacher, R.: Weak solutions of abstract evolutionary
integro-differential equations in Hilbert spaces. Funkcialaj
Ekvacioj {\bf 52} (2009), 1--18.

 }
\end{thebibliography}
\end{document}